\documentclass[english,british]{article}
\usepackage[T1]{fontenc}
\usepackage[latin9]{inputenc}
\usepackage{geometry}
\usepackage{color}
\usepackage{amsmath}
\usepackage{amssymb}

\makeatletter

\usepackage{amsthm}

\numberwithin{equation}{section}
\numberwithin{figure}{section}
\theoremstyle{plain}
\newtheorem{thm}{\protect\theoremname}[section]
 \theoremstyle{plain}
 \newtheorem{cor}[thm]{\protect\corollaryname}
 \newtheorem{remark}[thm]{Remark}

 \theoremstyle{definition}
 
 \newtheorem{example}[thm]{\protect\examplename}

\usepackage{babel}
\addto\captionsbritish{\renewcommand{\corollaryname}{Corollary}}
 \addto\captionsbritish{\renewcommand{\examplename}{Example}}
 \addto\captionsbritish{\renewcommand{\theoremname}{Theorem}}
 \addto\captionsenglish{\renewcommand{\corollaryname}{Corollary}}
 \addto\captionsenglish{\renewcommand{\examplename}{Example}}
 \addto\captionsenglish{\renewcommand{\theoremname}{Theorem}}
 \providecommand{\corollaryname}{Corollary}
 \providecommand{\examplename}{Example}
\providecommand{\theoremname}{Theorem}

\usepackage{babel}

\makeatother

\usepackage{babel}
\begin{document}

\title{A class of quasilinear second order partial differential equations
which describe spherical or pseudospherical surfaces}

\author{Diego Catalano Ferraioli$^{1}$, \quad{}Tarc\'{i}sio Castro Silva$^{2}$,
\quad{}Keti Tenenblat$^{3}$}
\date{}
\maketitle
\footnotetext[1]{Instituto de Matemática e Estatística- Universidade
Federal da Bahia, Campus de Ondina, Av. Adhemar de Barros, S/N, Ondina
- CEP 40.170.110 - Salvador, BA - Brazil, e-mail: diego.catalano@ufba.br.
Partially supported by CNPq, grant 310577/2015-2 and grant 422906/2016-6}
\footnotetext[2]{Department of Mathematics, Universidade de Brasilia,
Brasilia DF 70910-900, Brazil, e-mail: tarcisio@mat.unb.br. Partially
supported by FAPDF/Brazil, grant 0193.001346/2016, CAPES/Brazil-Finance Code 001 and by CNPq, grant
422906/2016-6}\footnotetext[3]{Department of Mathematics, Universidade
de Brasilia, Brasilia DF 70910-900, Brazil, e-mail: K.Tenenblat@mat.unb.br.
Partially supported by CNPq, grant 312462/2014-0, CAPES/Brazil-Finance Code 001 and FAPDF /Brazil,
grant 0193.001346/2016.} 

\begin{abstract}
Second order partial differential equations which describe spherical
surfaces (\textbf{ss}) or pseudospherical surfaces (\textbf{pss})
are considered. These equations are equivalent to the structure equations
of a metric with Gaussian curvature $K=1$ or $K=-1$, respectively,
and they can be seen as the compatibility condition of an associated
$\mathfrak{su}(2)$-valued or $\mathfrak{sl}\left(2,\mathbb{R}\right)$-valued
linear problem, also referred to as a zero curvature representation.
Under certain assumptions we give a complete and explicit classification
of equations of the form $z_{tt}=A(z,z_{x},z_{t})z_{xx}+B(z,z_{x},z_{t})z_{xt}+C(z,z_{x},z_{t})$
describing \textbf{pss} or \textbf{ss}, in terms of some arbitrary
differentiable functions. Several examples of such equations are provided
by choosing the arbitrary functions. In particular, well known equations
which describe pseudospherical surfaces, such as the short-pulse and
the constant astigmatism equations, as well as their generalizations and their spherical analogues
are included in the paper. 
\end{abstract}
\noindent 2010 \foreignlanguage{english}{\textit{Mathematics Subject
Classification}: 35G20, 53C21, 53B20}

\selectlanguage{english}%

\selectlanguage{british}%

\section{Introduction}

Equations describing spherical or pseudospherical surfaces are characterized
by the fact that their generic solutions provide metrics on a nonempty
open subsets of $\mathbb{R}^{2}$, with Gaussian curvature $K=1$
or $K=-1$, respectively. The best known example of equation describing
pseudospherical surfaces is the celebrated sine-Gordon equation (SG)
$z_{xt}=sin(z)$, which was discovered by Edmond Bour \cite{Baur}
to be equivalent to the Gauss-Codazzi equations for pseudospherical
surfaces in $\mathbb{R}^{3}$, in terms of Darboux asymptotic coordinates.
This equation is also the first well known example of an equation
integrable by means of geometric techniques originating from the classical
theory of surface transformation, and first applied to it by Bäcklund
and Bianchi. Nowadays SG is an important model in the description
of several nonlinear phenomena (see for example \cite{Josephson,Lamb,Skyrme}).

More in general the interest for equations describing spherical or
pseudospherical surfaces is motivated by their applications not only
in pure mathematics, but also in many problems of physics and applied
sciences. For instance, in the theory of solitons and integrable systems,
the interest on the general study of equations describing pseudospherical
surfaces began with the early observation made by Sasaki that ``all the
soliton equations in $1+1$ dimensions that can be solved by the AKNS
$2\times2$ inverse scattering method (for example, the sine-Gordon,
KdV or modified KdV equations) ... describe pseudospherical surfaces''
(see \cite{Sas1}). This motivated a general study of these equations,
initiated with the fundamental paper \cite{ChernTen} by S. S. Chern
and K. Tenenblat, which lead to important geometric interpretations
of Bäcklund transformations, conservation laws, non-local symmetries
and correspondences between solutions of different equations, for
equations of this class \cite{Beal-Rab-Ten,KT,Reyes-backlund,Rey9}.
The results of this study, together with the considerable effort addressed
over the past few decades to the possible applications of inverse
scattering method, gave a significant contribution to the discovery
of new integrable equations. For instance, Belinski-Zakharov system
in General Relativity \cite{Belinski}, the Jackiw-Teitelboim two-dimensional
gravity model \cite{Gegenberg-Kunstatter,Reyes-gravity}, the nonlinear
Schrödinger type systems \cite{Chen-Lee-Liu,DingTen}, the Rabelo's
cubic equation (also known as short-pulse equation) \cite{Beal-Rab-Ten,Rabelo,Sak1,Sak2,SW},
the Camassa-Holm, Degasperis-Procesi, Kaup\textendash{}Kupershmidt
and Sawada-Kotera equations \cite{CamHolm,Tarcisio-Keti,Rey6,Ray7,Rey8,Rey9}
are some important examples of equations describing pseudospherical
surfaces and integrable by inverse scattering method. On the other
hand, equations describing spherical surfaces are less studied, but
some known examples are the Schrödinger system $NLS^{+}$ (see \cite{DingTen}),
the Landau-Lifschitz equation for an isotropic chain (or isotropic
Heisenberg ferromagnet), Harry-Dym equation and Wadati-Konno-Ichikawa
equation (see \cite{WKI-2,WKI-1}).



In \cite{ChernTen} Chern and Tenenblat obtained characterization
results for evolution equations of the form $z_{t}=F(z,z_{x},...,\partial_{x}^{k}z)$,
$0\leq k$, under the assumptions that $f_{ij}=f_{ij}(z,z_{1},...,z_{k})$
and $f_{21}=\eta$, $\eta\in\mathbb{R}-\{0\}$. In the same paper
the authors also considered a similar problem for equations of the
form $z_{xt}=F(z,z_{x},...,\partial_{x}^{k}z)$. A noteworthy result
of this study was an effective method for the explicit determination
of entire new classes of differential equations that describe pseudospherical
surfaces. Motivated by the results of \cite{ChernTen}, in a series
of subsequent papers \cite{JorgeTen,Rabelo,RabTen,RabTen2}, the same
method was systematically implemented and new classes of pseudospherical
equations were identified still with the basic assumption that $f_{21}=\eta$.
Then in \cite{CavTen} the authors showed how the geometric properties
of pseudospherical surfaces may provide infinite number of conservation
laws when the functions $f_{ij}$ are analytic functions of the parameter
$\eta$. This parameter however is important not only for the existence
of infinite number of conservation laws, but is also related to the
existence of Bäcklund transformations and is crucial in the application
of inverse scattering method, as shown in \cite{Beal-Rab-Ten,ChernTen}.

In 1995, Kamran and Tenenblat \cite{KT} generalized some results
of \cite{ChernTen} by giving a complete characterization of evolution
equations of type $z_{t}=F(z,z_{x},...,\partial_{x}^{k}z)$ which
describe pseudospherical surfaces, in terms of necessary and sufficient
conditions that have to be satisfied by $F$ and the functions $f_{ij}=f_{ij}(z,z_{x},...,\partial_{x}^{k}z)$,
with no further additional conditions. Another generalization of some
aspects of \cite{ChernTen} came in 1998 by Reyes who considered in
\cite{Rey1} evolution equations of the more general form $z_{t}=F(x,t,z,z_{x},...,\partial_{x}^{k}z)$,
allowing $x,t$ to appear explicitly in the equation and assuming
that $f_{ij}=f_{ij}(x,t,z,z_{x},...,\partial_{x}^{k}z)$ and $f_{21}=\eta$.
Then, in a subsequent series of papers \cite{Rey6}-\cite{Rey8} Reyes
also studied other aspects of such equations.

In 2002, differential systems describing pseudospherical surfaces
or spherical surfaces (with constant positive curvature metrics) were
studied by Ding-Tenenblat in \cite{DingTen}. Such systems include
equations such as the nonlinear Schrödinger equation and the Heisenberg
Ferromagnet model, and large new families of differential systems
describing pseudospherical surfaces were obtained.

Also we mention that a higher dimensional geometric generalization
of the sine-Gordon equation, characterizing $n$-dimensional sub-manifolds
of the Euclidean $\mathbb{R}^{2n-1}$ with constant sectional curvature
$K=-1$, was considered in \cite{TT} and its intrinsic version as
a metric on a nonempty open subsets of $\mathbb{R}^{n}$, with $K=-1$,
was studied in \cite{BT}, by applying inverse scattering method.
Other differential $n$-dimensional systems that are the integrability
condition of linear systems of PDEs can be found in the so called
\textit{generating system} (see \cite{T} and its references).

The several characterization results obtained in \cite{ChernTen,KT,Rey1}
are extremely useful, either in checking if a given differential equation
describes pseudospherical surfaces or in generating large families
of such equations. For instance, as an application of \cite{KT},
Gomes \cite{Gomes} and Catalano-Tenenblat \cite{Catalano-Tenenblat}
classified evolution equations of the form $z_{t}=z_{xxxxx}+G(z,z_{x},z_{xx},z_{xxx},z_{xxxx})$
and $z_{t}=z_{xxxx}+G(z,z_{x},z_{xx},z_{xxx})$, respectively, under
the auxiliary assumptions that $f_{21}$ and $f_{31}$ are linear
combinations of $f_{11}$. Analogously, under the same assumption,
Castro Silva and Tenenblat in \cite{Tarcisio-Keti} have given a classification
of third order equations of the form $z_{t}=z_{xxt}+\lambda zz_{xxx}+G(z,z_{x},z_{xx})$,
with $\lambda\in\mathbb{R}$. Moreover in \cite{Catalano-Silva} it
has been given a classification of equations describing pseudospherical
surfaces of the form $z_{t}=A(x,t,z)z_{xx}+B(x,t,z,z_{x})$ with $f_{21}=\eta$,
$\eta\in\mathbb{R}$.

The results of \cite{Tarcisio-Keti,Catalano-Silva,Catalano-Tenenblat,Gomes}
permit the explicit description of huge classes of equations describing
pseudospherical surfaces which, apart from the already known examples,
represent a great amount of new equations whose physical relevance
is highly expected. For example, some applications of equations classified
by Rabelo and Tenenblat \cite{Beal-Rab-Ten,JorgeTen,Rabelo,RabTen}
have been recently discussed by Sakovich in a series of papers (see
for instance \cite{Sak2}). Of course the same should occur in the
case of results obtained in \cite{Tarcisio-Keti,Catalano-Silva,Catalano-Tenenblat,Gomes}
and the present paper.

In this paper, we classify second order equations (of parabolic, hyperbolic
or elliptic type) 
\begin{eqnarray}
z_{tt}=A(z,z_{x},z_{t})z_{xx}+B(z,z_{x},z_{t})z_{xt}+C(z,z_{x},z_{t}),\label{DT}
\end{eqnarray}
describing pseudospherical surfaces ($\delta=1$), or spherical surfaces
($\delta=-1$) with associated 1-forms $\omega_{i}=f_{i1}dx+f_{i2}dt$,
$1\leq i\leq3$, such that $f_{ij}$ are real differentiable functions
of $z_{t}$, $z$, $z_{x}$, $z_{xx}$, $\ldots$, $\partial_{x}^{k}z$,
satisfying 
\begin{eqnarray}
f_{21,z_{t}}=f_{31,z_{t}}=0.\label{condition}
\end{eqnarray}
The results of this classification are collected in Theorem \ref{teo3.2},
which provides an explicit description of such differential equations.
In particular, we point out special classes of equations in Corollaries
\ref{teo3.3}, \ref{teo3.4} and \ref{teo3.5}, which provide most
of the given explicit examples.

The paper is organized as follows. In Section \ref{sec:Preliminaries},
we collect some preliminaries on differential equations that describe
pseudospherical or spherical surfaces. Moreover, given such an equation
with associated $1$-forms $\omega_{i}$, we provide linear systems
of PDEs whose integrability condition is the differential equation.
Moreover, in Section \ref{section-mainresults}, we start providing some explicit examples which include generalizations of the constant astigmatism equation and of the short-pulse equation and then we state our main results, Theorem \ref{teo3.2} and  Corollaries \ref{teo3.3}-\ref{teo3.5}. In Section \ref{sec:Proofs}
we give the proof of Theorem \ref{teo3.2} and in Section \ref{further_ex}
we provide additional examples which follow from Theorem \ref{teo3.2}
and its Corollaries.

\section{\label{sec:Preliminaries}Preliminaries}

If $\left(M,\, g\right)$ is a 2-dimensional Riemannian manifold and
$\left\{ \omega_{1},\omega_{2}\right\} $ is a co-frame, dual to an
orthonormal frame $\left\{ e_{1},e_{2}\right\} $, then $g=\omega_{1}^{2}+\omega_{2}^{2}$
and $\omega_{i}$ satisfy the structure equations: $d\omega_{1}=\omega_{3}\wedge\omega_{2}$
and $d\omega_{2}=\omega_{1}\wedge\omega_{3}$, where $\omega_{3}$
denotes the connection form defined as $\omega_{3}(e_{i})=d\omega_{i}(e_{1},e_{2})$.
The Gaussian curvature of $M$ is the function $K$ such that $d\omega_{3}=-K\omega_{1}\wedge\omega_{2}$.

Now, a $k$-th order differential equation $\mathcal{E}$, for a scalar
or vector real-valued function $z\left(x,t\right)$ , \emph{describes
pseudospherical surfaces }\textbf{(pss)}\emph{, or spherical surfaces
}\textbf{(ss)}\emph{ }if it is equivalent to the structure equations
of a surface with Gaussian curvature $K=-\delta$, with $\delta=1$
or $\delta=-1$, respectively, i.e., 
\begin{equation}
\begin{array}{l}
d\omega_{1}=\omega_{3}\wedge\omega_{2},\quad d\omega_{2}=\omega_{1}\wedge\omega_{3},\quad d\omega_{3}=\delta\omega_{1}\wedge\omega_{2},\end{array}\label{struttura}
\end{equation}
where $\left\{ \omega_{1},\omega_{2},\omega_{3}\right\} $ are $1$-forms
\begin{equation}
\begin{array}{l}
\omega_{1}=f_{11}dx+f_{12}dt,\quad\omega_{2}=f_{21}dx+f_{22}dt,\quad\omega_{3}=f_{31}dx+f_{32}dt,\end{array}\label{eq:forms}
\end{equation}

\noindent such that $\omega_{1}\wedge\omega_{2}\neq0$, i.e., 
\begin{equation}
f_{11}f_{22}-f_{12}f_{21}\neq0,\label{eq:nondeg_cond}
\end{equation}
and $f_{ij}$ are functions of $x$, $t$, $z(x,t)$ and derivatives
of $z(x,t)$ with respect to $x$ and $t$.

Notice that, according to the definition, 
given a solution $z(x,t)$ of a 
\textbf{pss} (or \textbf{ss}) equation $\mathcal{E}$, we consider 
an open connected set   $U\subset\mathbb{R}^{2}$, contained in the domain of 
 $z(x,t)$, where the restriction
 of $\omega_{1}\wedge\omega_{2}$ to $z$ is everywhere nonzero on $U$.
  Such an open  set $U$ exists for generic solutions $z$.  
Thus, for generic solutions 
$z$ of a \textbf{pss} (or \textbf{ss}) $\mathcal{E}$, the restriction
$I[z]$ of $I=\omega_{1}^{2}+\omega_{2}^{2}$ to $z$ defines  a Riemannian metric $I[z]$, on $U$, with Gaussian
curvature $K=-\delta$. It is in this sense that one can say that
a \textbf{pss} (\textbf{ss}, resp.) describes, \textbf{pss} (\textbf{ss}, resp.). This is an intrinsic  geometric property of a  Riemannian metric 
(not immersed in an ambient space).

A classical example of equation describing \textbf{pss} is the sine-Gordon
equation $z_{xt}=sin\left(z\right)$, which corresponds to 
\[
\begin{array}{l}
\omega_{1}=\frac{1}{\eta}sin\left(z\right)\, dt,\qquad\omega_{2}=\eta\, dx+\frac{1}{\eta}cos\left(z\right)\, dt,\qquad\omega_{3}=z_{x}\, dx,\end{array}
\]
with $\eta\in\mathbb{R}-\{0\}$ (see \cite{ChernTen}). Observe that the 
solution $z_0(x,t)\equiv 0$ is not generic, in the sence that there is no open set $U$, where $\omega_{1}\wedge\omega_{2}\neq 0$, since $\omega_1(z_0)\equiv 0$.

Another classical example of equation describing \textbf{pss} is the
KdV equation $z_{t}=z_{xxx}+6zz_{x}$, which corresponds to 
\[
\begin{array}{l}
\omega_{1}=\left(1-z\right)\, dx+\left(-z_{xx}+\eta z_{x}-\eta^{2}z-2z^{2}+\eta^{2}+2z\right)\, dt,\vspace{4pt}\\
\omega_{2}=\eta\, dx+\left(\eta^{3}+2\eta z-2z_{x}\right)\, dt,\vspace{4pt}\\
\omega_{3}=-\left(1+z\right)\, dx+\left(-z_{xx}+\eta z_{x}-\eta^{2}z-2z^{2}-\eta^{2}-2z\right)\, dt,
\end{array}
\]
with $\eta\in\mathbb{R}$ (see \cite{ChernTen}).

Also the nonlinear dispersive wave equation (Camassa-Holm) 
\[
z_{t}-z_{xxt}=zz_{xxx}+2z_{x}z_{xx}-3zz_{x}-mz_{x},
\]
describes \textbf{pss} with 
\[
\begin{array}{l}
\omega_{1}=\left(z-z_{xx}+\frac{m+\eta^{2}}{2}-1\right)\, dx+\left[-z\left(f_{11}+1\right)\pm\eta z_{x}-\frac{m+\eta^{2}}{2}+1\right]\, dt,\vspace{4pt}\\
\omega_{2}=\eta\, dx+\left(-\eta z\pm z_{x}-\eta\right)\, dt,\vspace{4pt}\\
\omega_{3}=\pm\left(z-z_{xx}+\frac{m+\eta^{2}}{2}\right)dx+\left[\mp z\left(z-z_{xx}+\frac{m+\eta^{2}}{2}\right)+\eta z_{x}\mp z\mp\frac{m+\eta^{2}}{2}\right]\, dt,
\end{array}
\]
and $\eta\in\mathbb{R}$ (see \cite{Tarcisio-Keti}).

A classical example of equation describing \textbf{ss} is the nonlinear
Schrödinger system $NLS^{+}$ 
\[
\left\{ \begin{array}{l}
u_{t}+v_{xx}+2\left(u^{2}+v^{2}\right)u=0,\vspace{4pt}\\
-v_{t}+u_{xx}+2\left(u^{2}+v^{2}\right)v=0,
\end{array}\right.
\]
for the vector-valued function $\mathbf{z}(x,t)=(u(x,t),v(x,t))$.
This equation corresponds to 
\[
\begin{array}{l}
\omega_{1}=2vdx+\left(-4\eta v+2u_{x}\right)dt,\vspace{4pt}\\
\omega_{2}=2\eta dx+\left(-4\eta^{2}+2\left(u^{2}+v^{2}\right)\right)dt,\vspace{4pt}\\
\omega_{3}=-2udx+\left(2\eta u+2v_{x}\right)dt,
\end{array}
\]
with $\eta\in\mathbb{R}$ (see \cite{DingTen}).

Another example of equation describing \textbf{ss} is the Landau-Lifschitz
equation for an isotropic chain (or isotropic Heisenberg ferromagnet)
\[
\mathbf{z}_{t}=\mathbf{z}\wedge\mathbf{z}_{xx}
\]
for the vector-valued function $\mathbf{z}(x,t)=(u(x,t),v(x,t),w(x,t))$,
with $u^{2}+v^{2}+w^{2}=1$ and $\wedge$ denoting the cross product
in $\mathbb{R}^{3}$. This equation corresponds to 
\[
\begin{array}{l}
\omega_{1}=-2\eta vdx+\left(4\eta^{2}v+2\eta uw_{x}-2\eta wu_{x}\right)dt,\vspace{4pt}\\
\omega_{2}=-2\eta udx+\left(4\eta^{2}u-2\eta vw_{x}+2\eta wv_{x}\right)dt,\vspace{4pt}\\
\omega_{3}=2\eta wdx+\left(-4\eta^{2}w-2\eta vu_{x}+2\eta uv_{x}\right)dt,
\end{array}
\]
with $\eta\in\mathbb{R}-\{0\}$ (see \cite{DingTen}).

Equations which describe \textbf{pss}, or \textbf{ss}, can also be
characterized in few alternative ways. For instance, the system of
equations (\ref{struttura}) is equivalent to the integrability condition
\begin{equation}
d\Omega-\Omega\wedge\Omega=0,\label{eq:ZCR_omega}
\end{equation}
of the linear system 
\begin{equation}
dV=\Omega V,\label{struttura-2}
\end{equation}
for an auxiliary differentiable function $V=(v^{1},v^{2})^{T}$, with
$v^{i}=v^{i}\left(x,t\right)$, where $\Omega$ is either the $\mathfrak{sl}\left(2,\mathbb{R}\right)$-valued
1-form 
\[
\Omega=\frac{1}{2}\left(\begin{array}{cc}
\omega_{2} & \omega_{1}-\omega_{3}\\
\omega_{1}+\omega_{3} & -\omega_{2}
\end{array}\right),\qquad\text{when \;}\delta=1,
\]
or the $\mathfrak{su}\left(2\right)$-valued 1-form 
\[
\Omega=\frac{1}{2}\left(\begin{array}{cc}
i\omega_{2} & \omega_{1}+i\omega_{3}\\
-\omega_{1}+i\omega_{3} & -i\omega_{2}
\end{array}\right),\qquad\text{when \;}\delta=-1.
\]

\noindent Hence, for any solution $z=z\left(x,t\right)$ of $\mathcal{E}$,
defined on a domain $U\subset\mathbb{R}^{2}$, $\Omega$ is a Maurer-Cartan
form defining a flat connection on a trivial principal $SL\left(2,\mathbb{R}\right)$-bundle,
or $SU\left(2\right)$-bundle, over $\mathcal{U}$ (see for instance
\cite{DUB-KOM,SHARPE}).

Moreover, by using the matrices $X$ and $T$ such that $\Omega=Xdx+Tdt$
and $V:=\left(v^{1},v^{2}\right)^{T}$, (\ref{struttura-2}) can be
written as the linear problem 
\begin{equation}
\frac{\partial V}{\partial x}=XV,\quad\frac{\partial V}{\partial t}=TV.\label{eq:Prob_lin}
\end{equation}

\noindent It is easy to show that equations (\ref{struttura}) (or
(\ref{eq:ZCR_omega})) are equivalent to the integrability condition
of (\ref{eq:Prob_lin}), namely 
\begin{equation}
D_{t}X-D_{x}T+\left[X,T\right]=0,\label{eq:ZCR}
\end{equation}
where $D_{t}$ and $D_{x}$ are the total derivative operators with
respect to $t$ and $x$, respectively.

In the literature \cite{CraPirRob} 1-form $\Omega$, and sometimes
the pair $\left(X,T\right)$ or even (\ref{eq:ZCR}), is referred
to as an $\mathfrak{sl}\left(2,\mathbb{R}\right)$-valued, or $\mathfrak{su}\left(2\right)$-valued\emph{,
zero-curvature representation }for the equation $\mathcal{E}$. Moreover,
the linear system (\ref{struttura-2}) or (\ref{eq:Prob_lin}) is
usually referred to as \emph{the linear problem associated to $\mathcal{E}$. }
 It is this linear problem that, in some cases, is used in the construction
of explicit solutions of equations describing \textbf{pss}, or \textbf{ss},
by means of inverse scattering method \cite{AKNS,Beals,Beal-Coif,Beal-Rab-Ten,G.G.Krus.Miura.}.

We also notice that an equation which describes \textbf{pss} or \textbf{ss}
can alternatively be seen as the integrability condition 
\begin{equation}\label{eq:ZCR_omegahat}
d\hat{\Omega}-\hat{\Omega}\wedge\hat{\Omega}=0,
\end{equation}
of a linear system
\begin{equation}\label{num}  
d\psi =\hat{\Omega}\psi,
\end{equation}
for an auxiliary differentiable function $\psi=(\psi^1,\,\psi^2,\, \psi^3)^T$,  
with $\psi^i=\psi^i(x,t)$,  where   
\begin{eqnarray*}
\hat{\Omega}=\left(\begin{array}{ccc}
0 & \omega_{1} & \omega_{2}\\
\delta\omega_{1} & 0 & \omega_{3}\\
\delta\omega_{2} & -\omega_{3} & 0
\end{array}\right), \qquad \delta=\pm 1. \label{compeq3x3}
\end{eqnarray*}
By using the $3\times 3$ matrices $\hat{X}$ and $\hat{T}$ such that 
$\,\hat{\Omega}=\hat{X}dx+\hat{T}dt$,  then \eqref{num} can be written as the linear problem 
\begin{equation}
\label{linearhat} \frac{\partial\psi}{\partial x}=\hat{X}\psi,\qquad 
\frac{\partial\psi}{\partial t}=\hat{T}\psi.
\end{equation}
Notice that $\hat{X}$ and $\hat{T}$ are skew-symmetric when $\delta=-1$.
One can also easily see that \eqref{struttura}, and hence \eqref{eq:ZCR_omegahat}, are equivalent to the integrability condition 
of (\ref{linearhat}), i.e. $D_t\hat{X}-D_x\hat{T}+[\hat{X},\hat{T}]=0$.

\section{Some examples and main results \label{section-mainresults}}

In this section we will first illustrate some examples of equations
describing \textbf{pss} or \textbf{ss} of the form (\ref{eq:basic equation})
studied in the paper. Then, we will present the main classification
result, Theorem \ref{teo3.2}, and its useful Corollaries \ref{teo3.3}-\ref{teo3.5}.
The complete proof of Theorem \ref{teo3.2} is postponed to Section
\ref{sec:Proofs}, whereas further examples of equations described
by our main results are given in Section \ref{further_ex}.\\

\subsection{Examples \label{first_ex}}

In this subsection, we provide several examples.  Among them,  some
are already known in the literature. Further new examples will be
presented in Section \ref{further_ex}.

\medskip{}

\begin{example} \label{exe2.2} The following differential equation  
generalizes the constant astigmatism equation. 
\begin{equation}\label{gCAn}
z_{tt}=\frac{z_{xx}}{z^2}-\frac{(4z-3m)}{2(z-m)z^3}z_x^2+\frac{m}{2(z-m)z}z_t^2 -2+\frac{2m}{z}, 
\end{equation}
where $m\in\mathbb{R}$. It  is of the form (\ref{DT})
and it describes \textbf{pss} or \textbf{ss} (for $\delta=1$ or $\delta=-1$,
respectively) with 
\begin{equation}\label{eq:fij_genCAeqn}
\begin{array}{lll}
f_{11}=\displaystyle \frac{\eta\, z z_{t}+ \sqrt{\eta^2+\delta}\, z_x}{2\sqrt{z(z-m)}}, & f_{21}=\sqrt{z(\eta^2+\delta)},& f_{31}=\sqrt{z-m}, \vspace{6pt}\\
f_{12}=\displaystyle \frac{\eta z_x+\sqrt{\eta^2+\delta}\, z z_t}{2\sqrt{z(z-m)}\, z},& f_{22}=\displaystyle \frac{\eta}{\sqrt{z}},& f_{32}=0,
\end{array}
\end{equation}
where $\eta\in\mathbb{R}-\{0\}$. 
In particular, by taking $m=0$, equation (\ref{gCAn}) reduces
to the so called {\em constant astigmatism equation}  
\begin{equation}
z_{tt}=\frac{z_{xx}}{z^{2}}-2\frac{z_{x}^{2}}{z^{3}}-2,\label{CA_marvan}
\end{equation}
which was first studied in \cite{Baran-Marvan} by using the
alternative linear problem defined by 
\begin{equation}
\begin{array}{l}
f_{11}=\sqrt{\lambda^{2}+\lambda}\, z_{t},\vspace{5pt}\\
f_{21}=(\lambda+1)\, z^{-\lambda}+\lambda\, z^{\lambda+1},\vspace{5pt}\\
f_{31}=-\lambda\, z^{\lambda+1}+(\lambda+1)\, z^{-\lambda},
\end{array}\qquad\begin{array}{l}
f_{12}={\displaystyle \frac{\sqrt{\lambda^{2}+\lambda}}{z^{2}}}\, z_{x},\vspace{5pt}\\
f_{22}=\sqrt{\lambda^{2}+\lambda}\, z^{-\lambda-1}+\sqrt{\lambda^{2}+\lambda}\, z^{\lambda},\vspace{5pt}\\
f_{32}=\sqrt{\lambda^{2}+\lambda}\, z^{-\lambda-1}-\sqrt{\lambda^{2}+\lambda}\, z^{\lambda}.
\end{array}\label{eq:fij_CA_marvan}
\end{equation}

 \end{example}

\medskip{}

\begin{example} \label{exe2.1} The following differential equation  
generalizes the short-pulse equation. 
\begin{equation}
z_{tt}=\frac{2\,\delta}{z^{2}+m}\, z_{xt}-\frac{2z}{z^{2}+m}(z_{t}^{2}+1),\label{gSP}
\end{equation}
where $m\in\mathbb{R}$ and $\delta=\pm1$. It is of the form (\ref{DT})
and it describes \textbf{pss} or \textbf{ss} (for $\delta=1$ or $\delta=-1$,
respectively) with
\begin{equation}
\begin{array}{lll}
f_{11}={\displaystyle \frac{\lambda}{2}\delta\left(z^{2}+m\right)}z_{t},\hspace{25pt} & f_{21}={\displaystyle \frac{\delta\lambda}{2}\left(z^{2}+m\right)}+{\displaystyle \frac{1}{\lambda}},\qquad & f_{31}=\delta z,\vspace{5pt}\\
f_{12}={\displaystyle \lambda}z_{t}, & f_{22}=\lambda, & f_{32}=0,
\end{array}\label{fij_g_sp}
\end{equation}
where $\lambda\in\mathbb{R}-\{0\}$. In particular, in the \textbf{pss}
case ($\delta=1$), by taking $m=0$, equation (\ref{gSP}) reduces
to the so called {\em short-pulse equation} \cite{Beal-Rab-Ten} 
\begin{equation}
z_{tt}=2\frac{z_{xt}}{z^{2}}-\frac{2}{z}(z_{t}^{2}+1).\label{SP}
\end{equation}
The equation for $\delta=-1$ and $m=0$ was discussed in \cite{Sak1}.
\par \end{example}

\medskip{}

\begin{example}\label{ex_3.4.1} The differential equation 
\[
z_{tt}=m^2 z_{xx}+mpz^{p-1}z_x-pz^{p-1}z_t,\qquad p\in\mathbb{Z},
\]
where $m\in\mathbb{R}-\{0\}$, is of the form (\ref{DT}) and it describes
\textbf{pss} or \textbf{ss} (for $\delta=1$ or $\delta=-1$, respectively). 
It can also be written as 
\begin{equation*}
z_{tt}=m^{2}z_{xx}+m (z^p)_{x}- (z^p)_{t}.\label{eq_ex3.4.1}
\end{equation*}
The corresponding functions $f_{ij}$ are  
\begin{equation*}
\begin{array}{lll}
f_{11}=\lambda z_{t}+\lambda\left(m-\dfrac{1}{\lambda^{2}z^{2}-\delta}\right)z_{x}+\lambda\, z^{p}, & f_{21}=\dfrac{\eta}{\sqrt{|\lambda^{2}z^{2}-\delta|}},\hspace{10pt} & f_{31}=\lambda zf_{21}, \vspace{5pt}\\ 
f_{12}=\lambda\left(m-\dfrac{1}{\lambda^{2}z^{2}-\delta}\right)z_{t}+\lambda m^{2}\, z_{x}+\lambda\, mz^{p},\hspace{10pt} & f_{22}=mf_{21},& f_{32}=\lambda m zf_{21}, 
\end{array}
\end{equation*}
where $\eta,\lambda\in \mathbb{R}-\{0\}$.

\end{example}

\medskip{}

\begin{example}\label{ex_3.4.2n} The differential equation 
\begin{equation}\label{eq_ex3.4.2n}
z_{tt}=z^4z_{xx}+2z^3z_x^2+2\frac{z_t^2}{z}-\delta m^2 z
\end{equation}
is of the form (\ref{DT}) and describes   
\textbf{pss} or \textbf{ss} (for $\delta=1$ or $\delta=-1$, respectively)
 with  
\begin{equation}
\begin{array}{lll}
f_{11}=\displaystyle\frac{\delta}{m}( \frac{1}{z^2}z_t+z_x),\qquad & f_{21}=0,\qquad  & f_{31}=\displaystyle\frac{1}{z},\vspace{10pt}\\
f_{12}=\displaystyle \frac{\delta}{m}(z_t+z^2 z_x) & f_{22}=m, & f_{32}=z,
\end{array}\label{fij_eq_ex3.4.2n2}
\end{equation}
where $m\in \mathbb{R}-\{0\}$.
\end{example}

\medskip{}

\begin{example}\label{ex_3.4.2n2} The differential equation 
\begin{equation}\label{eq_ex3.4.2n2}
z_{tt}=z z_{xt}-z_x z_t+\frac{2}{z} z_t^2-\delta m^2z,
\end{equation}
is of the form (\ref{DT}) and describes
 \textbf{pss} or \textbf{ss} (for $\delta=1$ or $\delta=-1$, respectively)
 with 
\begin{equation}
\begin{array}{lll}
f_{11}=\displaystyle\frac{\delta}{mz^2}z_t,\qquad  & f_{21}=0,\qquad  & f_{31}=\frac{1}{z},\vspace{10pt}\\
f_{12}=\displaystyle\frac{\delta}{m z} z_t,\qquad  & f_{22}=m, & f_{32}=1,
\end{array}\label{fij_eq_ex3.4.2n}
\end{equation}
where $m\in \mathbb{R}-\{0\}$.
\end{example}

\medskip{}

\subsection{Main theorem  and corollaries}

We now state the main result of the paper

\begin{thm} \label{teo3.2} A second order differential equation
\begin{equation}
z_{tt}=A(z,z_{x},z_{t})z_{xx}+B(z,z_{x},z_{t})z_{xt}+C(z,z_{x},z_{t}),\label{eq:basic equation}
\end{equation}
describes \textbf{pss} or \textbf{ss} ($\delta=1$ or $\delta=-1$,
respectively), with correponding functions $f_{ij}=f_{ij}(z,z_{t},z_{x})$
such that $f_{21,z_{t}}=f_{31,z_{t}}=0$ if, and only if, 
\begin{equation}
A=\frac{f_{12,z_{x}}}{f_{11,z_{t}}}\qquad B=\frac{-f_{11,z_{x}}+f_{12,z_{t}}}{f_{11,z_{t}}}\qquad C=\frac{-z_{t}f_{11,z}+z_{x}f_{12,z}+\Delta}{f_{11,z_{t}}},\label{eq:eq_explicit}
\end{equation}
where 
\[
f_{11,z_{t}}\neq0,\qquad\qquad\Delta=f_{32}f_{21}-f_{31}f_{22},
\]
and in addition, depending on whether $\Delta\neq0$ or $\Delta\equiv0$,
the functions $f_{ij}$ have the form given below. \begin{description}\par 

\item {(a)} When $\Delta\neq0$, on a nonempty open set, then\par 
\begin{eqnarray*}
\left(\begin{array}{ccc}
f_{11}\\
\\
f_{12}
\end{array}\right)=\dfrac{1}{\Delta}\left(\begin{array}{ccc}
-\phi z_{x}-\psi_{21} &  & \varphi z_{x}+\psi_{31}\\
\\
-\phi z_{t}-\psi_{22} &  & \varphi z_{t}+\psi_{32}
\end{array}\right)\left(\begin{array}{ccc}
\psi_{21,z}\qquad-\psi_{22,z}\\
\\
\delta\psi_{31,z}\qquad-\delta\psi_{32,z}
\end{array}\right)\left(\begin{array}{ccc}
z_{t}\\
\\
z_{x}
\end{array}\right)
\end{eqnarray*}
where 
\[
\Delta=(\varphi\psi_{21}-\phi\psi_{31})z_{t}+(\phi\psi_{32}-\varphi\psi_{22})z_{x}+\psi_{32}\psi_{21}-\psi_{22}\psi_{31},
\]
\par 
\[
\begin{array}{l}
f_{21}=\phi z_{x}+\psi_{21},\vspace{5pt}\\
f_{31}=\varphi z_{x}+\psi_{31},
\end{array}\qquad\begin{array}{l}
f_{22}=\phi z_{t}+\psi_{22},\vspace{5pt}\\
f_{32}=\varphi z_{t}+\psi_{32},
\end{array}
\]
with $\phi$, $\varphi$ and $\psi_{ij}$, $i=2,3$, $j=1,2$, differentiable
functions of $z$ such that 
\[
(\phi z_{x}+\psi_{21})f_{12}-f_{11}(\phi z_{t}+\psi_{22})\neq0.
\]
\par 

\item {(b)} When $\Delta\equiv0$, then 
\begin{equation}
\begin{array}{lll}
f_{11}=h(z,z_{x},z_{t}),\vspace{5pt} & \quad\quad & f_{12}=\dfrac{1}{\phi z_{x}+\psi}\left[(\phi z_{t}+\chi)h+\delta z_{t}\left(\rho\psi\right)_{,z}-\delta z_{x}\left(\rho\chi\right)_{,z}\right],\vspace{5pt}\\
f_{21}=\phi z_{x}+\psi,\vspace{5pt} &  & f_{22}=\phi z_{t}+\chi,\vspace{5pt}\\
f_{31}=\rho f_{21}, &  & f_{32}=\rho f_{22},
\end{array}\label{eq:fij_Deltanulo}
\end{equation}
where $h_{z_{t}}\neq0$, $\phi$, $\psi$, $\chi$, $\rho$ are differentiable
functions of $z$ and one of the following two sub-cases occur: \begin{description}\par 

\item {(b.1)} when $\rho^{2}-\delta\neq0$, on a nonempty open
set, then 
\[
\psi=\frac{c_{1}}{\sqrt{\left|\rho^{2}-\delta\right|}},\qquad\chi=\frac{c_{2}}{\sqrt{\left|\rho^{2}-\delta\right|}},
\]
with $c_{1},\, c_{2}\in\mathbb{R}$, and $\rho_{,z}\neq0$, $c_{1}^{2}+c_{2}^{2}\neq0$
and $\phi^{2}+c_{1}^{2}\neq0$;\par 

\item {(b.2)} when $\rho^{2}-\delta\equiv0$, then $\delta=1$, $\rho=\pm1$,
$\left(\psi_{,z}\right)^{2}+\left(\chi_{,z}\right)^{2}\neq0$ and
$\phi^{2}+\psi^{2}\neq0$. \end{description} \end{description} \end{thm}
\medskip{}

According to Theorem \ref{teo3.2}, an equation of the form (\ref{eq:basic equation})
describes \textbf{pss} or \textbf{ss}, whith $f_{ij}$ satisfying  $f_{21,z_{t}}=f_{31,z_{t}}=0$
if, and only if, the coefficients $A$, $B$ and $C$ have the form
(\ref{eq:eq_explicit}) with $f_{ij}$ explicitly given as in the theorem.
Accordingly, one can distinguish two main cases $(a)\{\Delta\neq0\}$
and $(b)\{\Delta\equiv0\}$. In particular, the case $\Delta\equiv0$
subdivides in two further subcases $(b.1)$ and $(b.2)$. Therefore,
once a given equation is recognized to belong to one of these cases,
Theorem \ref{teo3.2} explicitly gives also the associated functions
$f_{ij}$.

\vspace{0.2in}

As an immediate consequence of our main result, in Corollaries \ref{teo3.3}-\ref{teo3.5}
below, we will provide special classes of differential equations included
in Theorem \ref{teo3.2}. Most of the examples of this paper are produced
by using these corollaries.

First, by considering $\phi=\varphi=0$ in Theorem \ref{teo3.2} (a), 
we observe that $f_{11,z_t}\neq 0$ is equivalent to $H=\psi_{31}^2-\delta\psi_{21}^2$ not being constant. 
Moreover, in this case $\Delta\neq 0$ reduces to $\Delta_{0}:=\psi_{32}\psi_{21}-\psi_{31}\psi_{22}\neq0$. 
Then one gets the following

\begin{cor} \label{teo3.3} Let $\psi_{ij}(z)$, $i=2,3$, $j=2,3$,
be differentiable functions such that $\Delta_{0}=\psi_{32}\psi_{21}-\psi_{31}\psi_{22}\neq0$,
$H=\psi_{31}^{2}-\delta\psi_{21}^{2}$ is not constant and $\delta=\pm1$. 
Then there is a differential equation describing \textbf{pss} or \textbf{ss}
($\delta=1$ or $\delta=-1$, respectively), of the form 
\[
z_{tt}=A\, z_{xx}+B\, z_{xt}+C,
\]
where 
\[
A={\displaystyle \frac{1}{H_{,z}}}\,(\delta \psi_{22}^2-\psi_{32}^2)_{,z},\qquad  B={\displaystyle \frac{2}{H_{,z}}}\,(\psi_{31}\psi_{32}-\delta\psi_{21}\psi_{22})_{,z},
\]
\[
C={\displaystyle \frac{2\delta\Delta_{0}}{H_{,z}}}(-f_{11,z}z_{t}+f_{12,z}z_{x}+\Delta_{0}).
\]
The corresponding functions $f_{ij}$ are 
\begin{eqnarray*}
\left(\begin{array}{ccc}
f_{11}\\
\\
f_{12}
\end{array}\right)=\dfrac{1}{\Delta_{0}}\left(\begin{array}{ccc}
-\psi_{21} &  & \psi_{31}\\
\\
-\psi_{22} &  & \psi_{32}
\end{array}\right)\left(\begin{array}{ccc}
\psi_{21,z} &  & -\psi_{22,z}\\
\\
\delta\psi_{31,z} &  & -\delta\psi_{32,z}
\end{array}\right)\left(\begin{array}{ccc}
z_{t}\\
\\
z_{x}
\end{array}\right)
\end{eqnarray*}
and $\, f_{ij}=\psi_{ij},\,$ for $i=2,3$,\; $j=1,2$. \end{cor}

\vspace{0.2in}

By choosing $\phi=0$ in Theorem \ref{teo3.2} (b.1), we 
observe that in this case $\rho^2-\delta\neq 0$ is equivalent to $\rho$  not being constant. Then considering $c_1=\eta\neq 0$ and  $m=c_2/c_1$,   one gets the following

\begin{cor} \label{teo3.4} Let $h(z,z_{x},z_{t})$ and $\rho(z)$
be differentiable functions such that $h_{,z_{t}}\neq0$ and $\rho$
is not constant. Then for any real number  $m$, 
there is a differential equation describing \textbf{pss} or \textbf{ss}
($\delta=1$ or $\delta=-1$, respectively), given by 
\[
z_{tt}=A\, z_{xx}+B\, z_{xt}+C,
\]
where 
\[
\begin{array}{l}
A=\displaystyle \frac{m}{h_{,z_{t}}}\left(h_{,z_{x}}+\frac{\rho'}{\rho^{2}-\delta}\right),\qquad B=m-\frac{1}{h_{,z_{t}}}\left(h_{z_{x}}+\frac{\rho'}{\rho^{2}-\delta}\right),\vspace{10pt}\\
C=\displaystyle \frac{ mz_x-z_t}{h_{z_t}}\left[h_z+z_x\left(\frac{\rho'}
{\rho^2-\delta}\right)'\right].
\end{array}
\]
The corresponding functions $f_{ij}$ are 
\begin{equation}
\begin{array}{lcl}
f_{11}=h,\qquad & f_{21}={\displaystyle \frac{\eta}{\sqrt{|\rho^{2}-\delta|}}}, & \qquad f_{31}=\rho f_{21},\vspace{10pt}\\
f_{12}= m h-{\displaystyle \frac{\rho'}{\rho^{2}-\delta}}(z_{t}-m z_{x}),\qquad & f_{22}=mf_{21}, & \qquad f_{32}=m \rho  f_{21}
\end{array}\label{eqcor3.4}
\end{equation}
where $\eta\in \mathbb{R}-\{0\}$.
\end{cor}

\vspace{.1in}

Finally, by considering $\phi=0$ in Theorem \ref{teo3.2} (b.2),
and observing that in this case $\psi(z)\neq0$, then one gets the
following

\vspace{0.2in}

\begin{cor} \label{teo3.5} Let $h(z,z_{t},z_{x})$, $\psi(z)\neq0$
and $\chi(z)$ be differentiable functions such that $h_{,z_{t}}\neq0$
and $\psi_{z}^{2}+\chi_{z}^{2}\neq0$. Then there is a differential
equation describing \textbf{pss}, given by 
\[
z_{tt}=A\, z_{xx}+B\, z_{xt}+C,
\]
where 
\[
A=\frac{1}{\psi h_{,z_{t}}}(\chi h_{,z_{x}}\mp\chi'),\qquad B=\frac{\chi}{\psi}-\frac{1}{h_{,z_{t}}}\left(h_{,z_{x}}\mp\frac{\psi'}{\psi}\right),
\]
\[
C=\frac{1}{h_{,z_{t}}}\left\{ -z_{t}h_{z}+z_{x}\left[\left(\frac{\chi}{\psi}\right)'h+\frac{\chi}{\psi}h_{z}\pm z_{t}\left(\frac{\psi'}{\psi}\right)'\mp z_{x}\left(\frac{\chi'}{\psi}\right)'\right]\right\} .
\]
The corresponding functions $f_{ij}$ are 
\begin{equation}
\begin{array}{lcl}
f_{11}=h,\qquad & f_{21}=\psi,\qquad & f_{31}=\pm\psi,\qquad\vspace{10pt}\\
f_{12}=\displaystyle\frac{1}{\psi}(\chi h\pm z_{t}\psi'\mp z_{x}\chi'),\qquad & f_{22}=\chi,\qquad & f_{32}=\pm\chi.
\end{array}\label{eqcor3.5}
\end{equation}
\par \end{cor}

\bigskip{}

\begin{remark}\label{remark obt exemplos} {\rm We observe that all the 
examples discussed in Subsection \ref{first_ex} follow from the above
results. In fact:

{\em i)  Example \ref{exe2.2}}  follows  from Corollary \ref{teo3.3}. In particular, equation (\ref{gCAn})
corresponds to the choices 
\[
\psi_{21}=\sqrt{z(\eta^2+\delta)}, \qquad \quad\psi_{22}=\frac{\eta}{\sqrt{z}},\qquad\quad  \psi_{31}= \sqrt{z-m},\qquad\quad  \psi_{32}=0, 
\]
where $m\in\mathbb{R}$, $\delta=\pm1$ and $\eta\in\mathbb{R}-\{0\}$.
Analogously, by choosing $\delta=1$, in Corollary \ref{teo3.3} and 
\[
\begin{array}{ll}
\psi_{21}=\left(\lambda+1\right)z^{-\lambda}+\lambda z^{\lambda+1},\vspace{5pt} & \psi_{22}=\sqrt{\lambda^{2}+\lambda}\left(z^{\lambda}+z^{-\lambda-1}\right)\vspace{5pt}\\
\psi_{31}=\left(\lambda+1\right)z^{-\lambda}-\lambda z^{\lambda+1}, & \psi_{32}=\sqrt{\lambda^{2}+\lambda}\left(-z^{\lambda}+z^{-\lambda-1}\right),
\end{array}
\]
with $\lambda\in\mathbb{R}-\{0\}$, one gets equation (\ref{CA_marvan}) 
and the corresponding functions $f_{ij}$ as in \eqref{eq:fij_CA_marvan}.

{\em ii)  Example \ref{exe2.1}} follows from Corollary \ref{teo3.3}
with the choices 
\[
\begin{array}{c}
\begin{array}{ll}
\psi_{21}= {\displaystyle \frac{1}{2}\delta\lambda\left(z^{2}+m\right)+\frac{1}{\lambda}},\qquad & \psi_{22}=\lambda,\qquad\psi_{31}=\delta z,\qquad\psi_{32}=0,\end{array}\end{array}
\]
where $m\in\mathbb{R}$, $\delta=\pm1$ and $\lambda\in\mathbb{R}-\{0\}$.\\

{\em iii) Example \ref{ex_3.4.1}}   follows from Corollary
\ref{teo3.4} with the choices 
\[
\begin{array}{l}
h=\lambda z_{t}+\lambda\left(m-\dfrac{1}{\lambda^{2}z^{2}-\delta}\right)z_{x}+\lambda\, z^{p},\qquad\rho=\lambda z,\end{array}
\]
where $\lambda,m\in\mathbb{R}-\{0\}$.\\

{\em iv) Examples \ref{ex_3.4.2n} and \ref{ex_3.4.2n2}}  follow from Corollary
\ref{teo3.3} by choosing  $\psi_{21}=0$, $\psi_{22}=m$, $\psi_{31}=1/z$ and 
$\psi_{32}=z$ and $\psi_{32}=1$ respectively, 
where $m\in\mathbb{R}-\{0\}$.}\\
\end{remark}

\section{A characterization result and proof of the main theorem\label{sec:Proofs}}

In order to simplify notations, $z_{1}$, $z_{2}$, ... will denote
partial derivatives $z_{x}$, $z_{xx}$, ... of $z$ with respect
to $x$. Accordingly equation \eqref{DT} will be rewritten as 
\begin{eqnarray}
z_{tt}=A(z,z_{1},z_{t})z_{2}+B(z,z_{1},z_{t})z_{1,t}+C(z,z_{1},z_{t}),\quad A^{2}+B^{2}+C^{2}\neq0.\label{class}
\end{eqnarray}

The proof of Theorem \ref{teo3.2} is based on the following preliminary
characterization. 

\begin{thm}\label{teo3.1} Equation \eqref{class} describes \textbf{pss}
or \textbf{ss} ($\delta=1$ or $\delta=-1$, respectively), with corresponding
functions $f_{ij}=f_{ij}(z,z_{t},z_{1},z_{2},...,z_{m})$, such that
$m\in\mathbb{N}$ and $f_{21,z_{t}}=f_{31,z_{t}}=0$, if, and only
if, the following conditions are satisfied: 
\begin{description}
\item [{(i)}] functions $f_{ij}$ only depend on $(z,z_{1},z_{t})$, in
particular $f_{11,z_{t}}\neq0$ and 
\begin{equation}
\begin{array}{l}
f_{21}=\phi z_{x}+\psi_{21},\vspace{5pt}\\
f_{31}=\varphi z_{x}+\psi_{31},
\end{array}\qquad\begin{array}{l}
f_{22}=\phi z_{t}+\psi_{22},\vspace{5pt}\\
f_{32}=\varphi z_{t}+\psi_{32},
\end{array}\label{cond3}
\end{equation}
with $\phi$, $\varphi$ and $\psi_{ij}$, $i=2,3$, $j=1,2$, differentiable
functions of $z$;
\item [{(ii)}] the functions $f_{11}$, $f_{12}$, $\phi$, $\varphi$,
$\psi_{ij}$, $A$, $B$ and $C$ satisfy the following equations
\begin{equation}
\left\{ \begin{array}{l}
-f_{11,z_{t}}A+f_{12,z_{1}}=0,\vspace{10pt}\\
-f_{11,z_{t}}B-f_{11,z_{1}}+f_{12,z_{t}}=0,\vspace{10pt}\\
-f_{11,z_{t}}\, C-z_{t}\, f_{11,z}+z_{1}f_{12,z}+(\varphi\psi_{21}-\phi\psi_{31})z_{t}+\vspace{3pt}\\
\qquad\qquad+(\phi\psi_{32}-\varphi\psi_{22})z_{1}+\psi_{32}\psi_{21}-\psi_{22}\psi_{31}=0,\vspace{10pt}\\
-z_{t}\,\psi_{21,z}+z_{1}\psi_{22,z}-f_{11}(\varphi z_{t}+\psi_{32})+(\varphi z_{1}+\psi_{31})f_{12}=0,\vspace{10pt}\\
-z_{t}\,\psi_{31,z}+z_{1}\psi_{32,z}+\delta[-f_{11}(\phi z_{t}+\psi_{22})+(\phi z_{1}+\psi_{21})f_{12}]=0,
\end{array}\right.\label{cond1-2-E1-E2-E3}
\end{equation}
with 
\begin{eqnarray}
(\phi z_{1}+\psi_{21})f_{12}-f_{11}(\phi z_{t}+\psi_{22})\neq0.\label{metric}
\end{eqnarray}

\end{description}
\end{thm}

\medskip{}

\begin{proof} Assume that $f_{ij}=f_{ij}(z,z_{t},z_{1},z_{2},...,z_{m})$
are differentiable functions such that $f_{21,z_{t}}=f_{31,z_{t}}=0$.
In this case the structure equations \eqref{struttura} reduce to\par 
\begin{equation}
\left\{ \begin{array}{l}
-f_{11,z_{t}}\, z_{tt}+{\displaystyle \sum_{k=0}^{m}f_{12,z_{k}}\, z_{k+1}+\left(f_{12,z_{t}}-f_{11,z_{1}}\right)z_{t,1}-\sum_{k=2}^{m}f_{11,z_{k}}\, z_{t,k}}\vspace{5pt}\\
\qquad-f_{11,z}\, z_{t}+f_{32}f_{21}-f_{31}f_{22}=0,\vspace{15pt}\\
{\displaystyle \sum_{k=0}^{m}f_{22,z_{k}}\, z_{k+1}+\left(f_{22,z_{t}}-f_{21,z_{1}}\right)z_{t,1}-\sum_{k=2}^{m}f_{21,z_{k}}\, z_{t,k}}\vspace{5pt}\\
\qquad-f_{21,z}\, z_{t}-f_{32}f_{11}+f_{31}f_{12}=0,\vspace{15pt}\\
{\displaystyle \sum_{k=0}^{m}}f_{32,z_{k}}\, z_{k+1}+\left(f_{32,z_{t}}-f_{31,z_{1}}\right)z_{t,1}-{\displaystyle \sum_{k=2}^{m}}f_{31,z_{k}}\, z_{t,k}\vspace{5pt}\\
\qquad-f_{31,z}\, z_{t}+\delta\left(f_{21}f_{12}-f_{22}f_{11}\right)=0.
\end{array}\right.\label{eqstr_carat}
\end{equation}
Therefore equation \eqref{class} is equivalent to the structure equations
\eqref{struttura}, with $\omega_{1}\wedge\omega_{2}\neq0$, if, and
only if, it is equivalent to \eqref{eqstr_carat} with $f_{11}f_{22}-f_{12}f_{21}\neq0$
and in addition $f_{11,z_{t}}\neq0$. Then it is readily shown that
\eqref{eqstr_carat}, together with $f_{11}f_{22}-f_{12}f_{21}\neq0$,
is equivalent to \textit{(i-ii)}.\\
 Indeed, in \eqref{eqstr_carat}, by taking the coefficients of $z_{t,k}$,
with $k\geq2$, one gets that $f_{11}$, $f_{21}$ and $f_{31}$ do
not depend on $z_{k}$. Analogously, by taking the coefficients of
$z_{k+1}$, with $k\geq2$, one gets that $f_{12}$, $f_{22}$ and
$f_{32}$ do not depend on $z_{k}$. In particular, by taking the
coefficients of $z_{2}$ in the second and third equation of \eqref{eqstr_carat}
one also gets that $f_{22}$ and $f_{32}$ do not depend on $z_{1}$.\\
 On the other hand, by taking the coefficients of $z_{t,1}$ in the
second and third equation of \eqref{eqstr_carat} one also gets that\par 
\begin{eqnarray}
f_{21,z_{1}}-f_{22,z_{t}}=0, & f_{31,z_{1}}-f_{32,z_{t}}=0.\label{D2-a}
\end{eqnarray}
Thus $f_{21}$, $f_{31}$, $f_{22}$ and $f_{32}$ are as in \eqref{cond3},
and $f_{11}f_{22}-f_{12}f_{21}\neq0$ reduces to \eqref{metric}.\\
 Finally, \eqref{cond1-2-E1-E2-E3} follows by taking the coefficients
of $z_{t,1}$ and $z_{2}$ in the first equation of \eqref{eqstr_carat},
and further using \eqref{cond3} in \eqref{eqstr_carat}. The converse
is a straightforward computation.\par \end{proof}

\medskip{}

\subsection*{Proof of Theorem \ref{teo3.2}}

In view of Theorem \ref{teo3.1} (\textsl{i}), since $f_{11,z_{t}}\neq0$,
the first three equations of \eqref{cond1-2-E1-E2-E3} give us $A$,
$B$ and $C$ as in \eqref{eq:eq_explicit}, where $\Delta:=f_{32}f_{21}-f_{31}f_{22}$,
with $f_{21}$, $f_{31}$, $f_{22}$ and $f_{32}$ given by \eqref{cond3}.

Now, if $\Delta$ is nonzero, one can resolve the last two equations
of system \eqref{cond1-2-E1-E2-E3} in terms of $f_{11}$ and $f_{12}$.
Thus, we prove Theorem \ref{teo3.2} item (a). 

On the other hand, assuming that $\Delta\equiv0$, we can easily prove
that $f_{21}\neq0$. Indeed, if on the contrary $f_{21}=0$, then
by \eqref{cond3} one should have that $\phi=\psi_{21}=0$ and, in
view of \eqref{metric}, that $\psi_{22}\neq0$. Furthermore, $\Delta=0$,
$\phi=\psi_{21}=0$ and $\psi_{22}\neq0$ would imply that $\psi_{31}=0$
and from the last equation of \eqref{cond1-2-E1-E2-E3} one would get 
\begin{equation}
z_{1}\psi_{32,z}-\delta f_{11}\psi_{22}=0.\label{eq:aux1}
\end{equation}
Then, by differentiating (\ref{eq:aux1}) with respect to $z_{t}$,
one would get $f_{11,z_{t}}\psi_{22}=0$ and hence $f_{11,z_{t}}=0$,
which is a contradiction. Therefore, we can assume $f_{21}\neq0$
and, in view of \eqref{cond3}, we have $f_{22,z_{t}}=\phi=f_{21,z_{1}}$
and $f_{32,z_{t}}=\varphi=f_{31,z_{1}}$.

Now, by differentiating $\Delta=0$ with respect to $z_{t}$, one
has 
\begin{equation}
0=f_{31}f_{22,z_{t}}-f_{21}f_{32,z_{t}}=f_{31}f_{21,z_{1}}-f_{21}f_{31,z_{1}},\label{eq:aux2}
\end{equation}
hence $(f_{31}/f_{21})_{,z_{1}}=0$ and $f_{31}=\rho f_{21}$, with
$\rho$ a differentiable function of $z$.

It follows that, in view of (\ref{eq:aux2}), $f_{32}=\rho f_{22}$
and \eqref{cond3} entails that $\varphi=\rho\phi$ and $\psi_{3j}=\rho\psi_{2j}$,
$j=1,2$. In the statement of the main result, when $\Delta=0$, we
considered $\psi_{21}=\psi$ and $\psi_{22}=\chi$, in order to simplify
the notation.

Therefore, since $f_{21}\neq0$, by last equation of \eqref{cond1-2-E1-E2-E3}
one gets 
\begin{eqnarray}
f_{12}=\frac{1}{\phi z_{1}+\psi_{21}}[(\phi z_{t}+\psi_{22})f_{11}+\delta z_{t}(\rho\psi_{21})_{,z}-\delta z_{1}(\rho\psi_{22})_{,z}],\label{f12}
\end{eqnarray}
which substituted into the fourth equation of \eqref{cond1-2-E1-E2-E3}
leads to 
\begin{eqnarray}
z_{t}[-\psi_{21,z}+\delta\rho(\rho\psi_{21})_{,z}]+z_{1}[\psi_{22,z}-\delta\rho(\rho\psi_{22})_{,z}]=0,\label{sub_f12}
\end{eqnarray}
Now, by differentiating \eqref{sub_f12} with respect to $z_{t}$
and $z_{1}$, one gets respectively 
\begin{eqnarray*}
-\psi_{21,z}+\delta\rho(\rho\psi_{21})_{,z}=0,\quad\psi_{22,z}-\delta\rho(\rho\psi_{22})_{,z}=0,
\end{eqnarray*}
that is, 
\begin{eqnarray}
(\rho^{2}-\delta)\psi_{21,z}+\rho\,\rho_{,z}\psi_{21}=0,\quad(\rho^{2}-\delta)\psi_{22,z}+\rho\,\rho_{,z}\psi_{22}=0.\label{last_sis}
\end{eqnarray}

Now, if $\rho^{2}-\delta\neq0$ then \eqref{last_sis} implies that
 $\psi_{21}=c_{1}|\rho^{2}-\delta|^{-1/2}$ and $\psi_{22}=c_{2}|\rho^{2}-\delta|^{-1/2}$,
where $c_{1}$ and $c_{2}$ are real constants. Furthermore, the non-degeneracy
condition \eqref{metric} reads 
\[
z_{t}(\rho\psi_{21})_{,z}-z_{1}(\rho\psi_{22})_{,z}=(c_{1}z_{t}-c_{2}z_{1})(\rho|\rho^{2}-\delta|^{-1/2})_{,z}\neq0.
\]
This concludes the proof of Theorem \ref{teo3.2} item (b) subcase
\foreignlanguage{english}{\textbf{(b.1)}.}

Otherwise, if $\rho^{2}-\delta=0$ then 
$\delta=1$ and  $\rho=\pm1$. In this case, \eqref{last_sis} is trivially satisfied  
and the non-degeneracy condition \eqref{metric} reads $z_{t}\psi_{21,z}-z_{1}\psi_{22,z}\neq0$.
This concludes the proof of Theorem \ref{teo3.2} item (b) subcase
\foreignlanguage{english}{\textbf{(b.2)}.}

The converse is a straightforward computation in each case. \hfill\(\Box\)

\section{Further examples \label{further_ex}}

Here we provide some additional examples given by Theorem \ref{teo3.2},
and its Corollaries \ref{teo3.3}-\ref{teo3.5}.

\medskip{}
\begin{example} \label{ex_cor} The differential equation
\[
z_{tt}=m(m-\ell'(z))z_{xx}+\ell'(z)z_{xt}-m\ell''(z)z_x^2+\ell''(z)z_xz_t
\]
or equivalently
\begin{equation}\label{eq.3.3n}
z_{tt}=m^2z_{xx}-m(\ell(z))_{xx}+(\ell(z))_{xt}
\end{equation}
where $ m\in\mathbb{R}-\{0\}$ and $\ell(z)$ is a differentiable
function, can be obtained from Corollary \ref{teo3.4}, by choosing
\begin{eqnarray*}
h(z,z_{x},z_{t})=\lambda z_{t}+\left(m\lambda-\frac{\rho'}{\rho^{2}-\delta}-\lambda\ell'\right)z_{x},
\end{eqnarray*}
with $\lambda\in\mathbb{R}-\{0\}$, $\delta=\pm1$ and 
$\rho(z)$ is an arbitraty non constant differentiable function. With these
choices, (\ref{eq.3.3n}) describes \textbf{pss} or \textbf{ss} (for
$\delta=1$ or $\delta=-1$, respectively) with 
\begin{equation}
\begin{array}{lll}
f_{11}= h, \qquad & f_{21}={\displaystyle \frac{\eta}{\sqrt{|\rho^{2}-\delta|}}},\qquad & f_{31}=\rho f_{21},\vspace*{.3cm}\\
 f_{12}= \lambda m(m-\ell')z_x +(m\lambda-\displaystyle\frac{\rho'}{\rho^2-\delta})z_t,
 \qquad & f_{22}=m f_{21},\qquad & f_{32}=m\rho f_{21},
\end{array}
\end{equation}
where $\eta\in\mathbb{R}-\{0\}$.
\end{example} 
Notice that the  linear problem associated to \eqref{eq.3.3n} depends on a  parameter $\eta$ and  it also involves a non constant arbitrary 
function $\rho(z)$. 

\medskip{}

\begin{example} \label{ex_cor-1} The differential equation for $z(x,t)$  
\begin{equation}\label{eq.3.3-1}
z_{tt}=\alpha(z)z_{xx}+\ell(z)z_{xt}+\alpha'(z)z_x^2+\ell'(z)z_xz_t,
\end{equation}
or equivalently 
\begin{equation*}\label{eq.3.3-1n}
z_{tt}=(\alpha(z)z_x)_x+(\ell(z)z_x)_t, 
\end{equation*}
where $\alpha(z)$ and $\ell(z)$ are a differentiable functions, can be obtained from Corollary \ref{teo3.5}, by choosing 
\[
h(z,z_{x},z_{t})=\lambda z_{t}, \qquad  \lambda\in\mathbb{R}-\{0\},
\]
and $\psi(z)$, $\chi(z)$ satisfying $\psi\neq0$, $(\chi/\psi)'\neq0$ and the following linear system of ODE
\begin{eqnarray*}
\left(\begin{array}{ccc}
\chi'\\
\\
\psi'
\end{array}\right)=\mp \lambda\left(\begin{array}{ccc}
0 &  & \alpha\\
\\
1 &  & -\ell
\end{array}\right)\left(\begin{array}{ccc}
\chi\\
\\
\psi
\end{array}\right).
\end{eqnarray*}
With these choices, (\ref{eq.3.3-1}) describes \textbf{pss} with 
\begin{equation*}
\begin{array}{lcl}
f_{11}=\lambda\, z_{t},\qquad & f_{21}=\psi,\qquad & f_{31}=\pm\psi,\qquad\vspace{10pt}\\
f_{12}=\lambda\ell\, z_t + \lambda \alpha\, z_x,\qquad & f_{22}=\chi,\qquad & f_{32}=\pm\chi.
\end{array}
\end{equation*}
The linear problem associated to this equation depends not only on the parameter $\lambda$ but also on the parameters given by the initial conditions chosen for the  solution of  the 
system of ODEs for $\psi$ and $\chi$.
\end{example} \medskip{}

\remark {\rm By considering $\ell(z)=0$ in \eqref{eq.3.3-1}, the resulting
differential equation appers in \cite{EG}}.

\medskip{}

\begin{example} \label{ex_k2} The differential equation 
\begin{equation}
\begin{array}{l}
z_{tt}={\displaystyle -\frac{p}{q}z^{2(p-q)}z_{xx}+\left(\frac{p}{q}+1\right)z^{p-q}z_{xt}-\frac{p}{q}\left(2p-q-1\right)z^{2(p-q)-1}z_{x}^{2}+}\vspace{10pt}\\
\qquad\qquad\qquad\qquad\qquad{\displaystyle +\frac{(p-1)(p+q)}{q}z^{p-q-1}z_{x}z_{t}-(q-1)\frac{z_{t}^{2}}{z}+\frac{\delta m^{2}}{q}z}
\end{array}\label{eqk2}
\end{equation}
where $m\in\mathbb{R}-\{0\}$, $q\neq 0,p \in\mathbb{Z}$ can be obtained
from Corollary \ref{teo3.3}, by choosing 
\[
\psi_{21}=0,\qquad\psi_{22}=m,\qquad\psi_{31}=z^{q},\qquad\psi_{32}=z^{p}.
\]
With these choices, (\ref{eqk2}) describes\textbf{ pss} or \textbf{ss}
(for $\delta=1$ or $\delta=-1$, respectively) with 
\begin{equation*}
\begin{array}{lll}
f_{11}={\displaystyle \frac{\delta p}{m}z^{p-1}z_{x}-\frac{\delta q}{m}z^{q-1}z_{t}},\qquad & 
f_{21}=0,\qquad &
f_{31}=z^{q},\vspace*{0.3cm}\\
f_{12}={\displaystyle \frac{\delta p}{m}z^{2p-q-1}z_{x}-\frac{\delta q}{m}z^{p-1}z_{t}},\qquad & 
f_{22}=m,\qquad &
f_{32}=z^{p}.
\end{array}
\end{equation*}

Quite simple instances of (\ref{eqk2}) correspond to the choices:
\emph{(i)} $q=1$ and $p=-2,-1,0,1,2$; \emph{(ii)} $q=2p-1$ and $p=1,2$; or $q=-1$ and $p=0$ (see Example \ref{ex_3.4.2n2}); 
\emph{(iii)} $p=q$; \emph{(iv)} $p+q=0$ in particular $p=1=-q$ gives Example \ref{ex_3.4.2n}. Observe that choosing $p=q=1$ (\ref{eqk2}) reduces to a linear equation. \end{example} 

\medskip{}

\begin{example} \label{ex_k3} The differential equation 
\begin{equation}\label{eqk3n}
z_{tt}=(-\frac{1}{2}z^2+m)z_{xx}+(-\frac{3}{2}z+\frac{m}{z})z_{xt} 
+(-\frac{3}{2}+\frac{m}{z^2})z_tz_x- zz_x^2+(4\delta-1)(m_2-2m)^2z.
\end{equation}
where $m ,m_{2}\in\mathbb{R}$, with $m_{2}\neq 2m$, can be
obtained from Corollary \ref{teo3.3}, by choosing 
\[
\psi_{21}=z,\qquad\psi_{22}=-\frac{z^{2}}{2}+2\delta m_2-(4\delta-1)m,\qquad\psi_{31}=2z,\qquad\psi_{32}=-z^{2}+m_{2}.
\]
With these choices, (\ref{eqk3n}) describes\textbf{ pss} or \textbf{ss}
(for $\delta=1$ or $\delta=-1$, respectively) with 
\begin{equation*}
\begin{array}{lll}
f_{11}=\displaystyle \frac{-1}{(m_2-2m)}(zz_x+z_t)\qquad & 
f_{21}=z,\qquad & f_{31}=2z,\vspace*{.3cm}\\
f_{12}=\displaystyle (-\frac{z}{2}+\frac{m}{z} )f_{11}
\qquad & 
f_{22}=-\displaystyle\frac{z^{2}}{2}+2\delta m_2-(4\delta-1)m,\qquad & 
 f_{32}=-z^2+m_2. 
\end{array}
\end{equation*}

A quite simple instance of (\ref{eqk3n}) corresponds to the choice
$m_{2}=\left(4-\delta\right)m/2$, where $m\in\mathbb{R}$.\end{example}

\medskip{}

\begin{example} \label{ex_k4} The differential equation 
\begin{equation}
z_{tt}=m_{1}m_{2}z_{xx}+(m_{1}-m_{2})z_{xt}+\left(m_{1}z_{x}-z_{t}\right)\ell'(z),\label{eq4k}
\end{equation}
where $m_{1},m_{2}\in\mathbb{R}$ and $\ell=\ell(z)$ is a differentiable
function, can be obtained from Corollary \ref{teo3.5}, by choosing
\[
h=z_{t}+\left(m_{2}\pm{\displaystyle \frac{1}{z}}\right)z_{x}+\ell,\qquad\psi=\eta\, z,\qquad\chi=\eta\, m_{1}z,
\]
where $\eta\in\mathbb{R}-\{0\}$. With these choices, (\ref{eq4k})
describes\textbf{pss} 
with 
\[
\begin{array}{lll}
f_{11}=z_{t}+\left(m_{2}\pm{\displaystyle \frac{1}{z}}\right)z_{x}+\ell,\qquad &f_{21}=\eta\, z,\qquad & f_{31}=\pm\eta\, z,\vspace*{0.3cm}\\
f_{12}=m_{1}m_{2}z_{x}+\left(m_{1}\pm{\displaystyle \frac{1}{z}}\right)z_{t}+m_{1}\ell,\qquad &
f_{22}=\eta\, m_{1}z,\qquad&
f_{32}=\pm\eta\, m_{1}z.
\end{array}
\]
The linear sysytem associated to this equations depends on the parameter 
$\eta$.
\end{example} \medskip{}

\medskip{}
\begin{example} \label{ex_gen} The differential equation 
\begin{equation}
z_{tt}=-\frac{z_{t}}{zz_{x}}z_{xx}+\frac{z\, z_{t}+z_{x}}{z\, z_{x}}z_{xt}\mp\frac{\left(z\, z_{t}-z_{x}\right)^{3}}{z\, z_{x}}\label{eq_gen}
\end{equation}
can be obtained from Theorem \ref{teo3.2}, by choosing 
\[
\phi=1,\qquad\varphi=\psi_{21}=\psi_{22}=0,\qquad\psi_{31}=\eta\, z,\qquad\psi_{32}=\eta,
\]
where $\eta\in\mathbb{R}-\{0\}$. With these choices, (\ref{eq_gen})
describes\textbf{ ss} (for $\delta=1$ or $\delta=-1$, respectively)
with
\[
f_{11}={\displaystyle \frac{\delta\eta\, z\, z_{t}}{z_{x}-z\, z_{t}}},\qquad f_{12}={\displaystyle \frac{\delta\eta\, z_{t}}{z_{x}-z\, z_{t}}},\qquad f_{21}=z_{x},\qquad f_{22}=z_{t},\qquad f_{31}=\eta\, z,\qquad f_{32}=\eta.
\]
The linear  problem associated to this equation depends on the parameter $\eta$. 
 \end{example} \medskip{}

\end{document}